\documentclass[a4note,12pt,leqno]{article}

\usepackage{amssymb}

\begin{document}

\title{Approximation for extinction probability of the contact process based on the Gr\"obner basis}

\author{Norio Konno \\
Department of Applied Mathematics\\
Yokohama National University\\
\\}
\date{}
\maketitle

\newtheorem{thm}{Theorem}[section]
\newtheorem{lem}[thm]{Lemma}
\newtheorem{concl}[thm]{Claim}
\newtheorem{cor}[thm]{Corollary}
\newtheorem{pro}[thm]{Proposition}
\newtheorem{conj}[thm]{Conjecture}

\newcommand{\nl}{\nu_{\lambda}}
\newcommand{\zd}{{\bf Z}^d}
\newcommand{\rl}{\rho_{\lambda}}
\newcommand{\rtl}{\rho_{t,\lambda}}
\newcommand{\nlt}{\overline{\rho}_{\lambda,t}}
\newcommand{\lam}{\lambda}
\newcommand{\dt}{{\partial \over \partial t}}
\newcommand{\omu}{\overline{\mu}}
\newcommand{\ep}{\varepsilon}
\newcommand{\sn}{S^{(n)}}
\newcommand{\kuro}{\bullet}
\newcommand{\siro}{\circ}
\newcommand{\batu}{\times}
\newcommand{\ket}[1]{|#1\rangle}
\newcommand{\bra}[1]{\langle#1|}
\newcommand{\U}{\bar{U}}
\newcommand{\braa}{\langle}
\newcommand{\kett}{\rangle}
\newcommand{\RM}{\mathbb{R}}
\newcommand{\ZM}{\mathbb{Z}}
\newcommand{\ZMP}{\mathbb{Z}_{+}}
\newcommand{\QM}{\mathbb{Q}}
\newcommand{\NM}{\mathbb{N}}
\newcommand{\CM}{\mathbb{C}}
\newcommand{\TM}{\mathbb{T}}
\newcommand{\KM}{\mathbb{K}}

\par\noindent
\begin{small}
{\bf Abstract}. In this note we give a new method for getting a series of approximations for the extinction probability of the one-dimensional contact process by using the Gr\"obner basis.
\end{small}

\section{Introduction}
Let $X = \{ 0,1 \}^{\ZM^d}$ denote a configuration space, where $\ZM^d$ is the $d$-dimensional integer lattices. The contact process $\{ \eta_t : t \ge 0 \}$ is an $X$-valued continuous-time Markov process. The model was introduced by Harris in 1974 \cite{h1974} and is considered as a simple model for the spread of a disease with the infection rate $\lam$. In this setting, an individual at $x \in \ZM^d$ for a configuration $\eta \in X$ is infected if $\eta (x)=1$ and healthy if $\eta (x)=0$. The formal generator is given by
\begin{eqnarray*}
\Omega f(\eta)=\sum_{x \in \ZM^d} c(x,\eta)[f(\eta^{x})-f(\eta)],
\end{eqnarray*}
where $\eta^x \in X$ is defined by $\eta^x (y)=\eta (y) \> (y \ne x),$ and 
$\> \eta^x (x)=1- \eta(x)$. Here for each $x \in \ZM^d$ and $\eta \in X,$ the transition rate is 
\begin{eqnarray*} 
c(x,\eta)= (1- \eta (x)) \times \lambda \sum_{y:\mid y-x
\mid=1} \eta (y) + \eta (x),
\end{eqnarray*}
with $|x| = |x_1| + \cdots + |x_d|.$ In particular, the one-dimensional contact process is 
\begin{eqnarray*}
001 \to 011 && \qquad \qquad \hbox{at rate} \qquad  \lam,
\\
100 \to 110 && \qquad \qquad \hbox{at rate} \qquad  \lam,
\\
101 \to 111 && \qquad \qquad \hbox{at rate} \qquad  2 \lam,
\\
1 \> \> \to \> \> 0 \>  \> && \qquad \qquad \hbox{at rate}
\qquad 1. 
\end{eqnarray*}
Let $Y= \{ A \subset \ZM^d : |A| < \infty \}$, where $|A|$ is the number of elements in $A$. Let $\xi^A _t (\subset \ZM^d)$ denote the state at time $t$ of the contact process with $\xi^A _0 = A$. There is a one-to-one correspondence between $\xi^A _t (\subset \ZM^d)$ and $\eta_t \in X$ such that $x \in \xi^A_t$ if and only if $\eta_t (x) = 1$. For any $A \in Y$, we define the extinction probability of $A$ by $\lim_{t \to \infty} P(\xi^A _t = \emptyset).$ Define $\nl (A) = \nl \{ \eta \> : \> \eta (x)=0 \>$ for any $\> x \in A \}$, where $\nl$ is an invariant measure of the process starting from a configuration: $\eta (x) =1 \> (x \in \ZM^d)$ and is called the {\it upper invariant measure}. In other words, let $\delta_1 S(t)$ denote the probability measure at time $t$ for initial probability measure $\delta_i$ which is the pointmass $\eta \equiv i \> (i=0,1).$ Then $\nl = \lim_{t \to \infty} \delta_1 S(t)$. Then self-duality of the process implies that $\nl (A) = \lim_{t \to \infty} P(\xi^A _t = \emptyset)$. The correlation identities for $\nl (A)$ can be obtained as follows:
\begin{thm}
\label{thm1} 
For any $ A \in Y$, 
\begin{eqnarray*}
\lambda \sum_{x \in A} \sum_{y : \mid y-x \mid = 1} \Bigl[ \nl
(A \cup \{ y \}) - \nl (A) \Bigr] + \sum_{x \in A} \Bigl[ \nl (
A \setminus \{ x \}) - \nl (A) \Bigr] =0.
\end{eqnarray*}
\end{thm}
From now on we consider the one-dimensional case. We introduce the following notation:
\begin{eqnarray*}
\nl (\siro)=\nl (\{0\}), \> \nl (\siro \siro)=\nl (\{0,1\}),
\> \nl(\siro \batu \siro)=\nl (\{0,2\}), \> \ldots .
\end{eqnarray*}
By Theorem {\rmfamily \ref{thm1}}, we obtain 
\begin{cor}
\label{cor2}
\begin{eqnarray}
&& 
2 \lam \nl (\siro \siro) -(2 \lam+1) \nl (\siro) + 1=0,
\\
&&
\lam \nl (\siro \siro \siro) - (\lam +1)\nl (\siro \siro) + \nl (\siro) =0,
\\
&&
2 \lam \nl (\siro \siro \siro \siro) + \nl (\siro \batu \siro) - ( 2 \lam +3 )\nl (\siro \siro \siro) + 2 \nl (\siro \siro) =0,
\\
&&
\lam \nl (\siro \siro \batu \siro) - (2 \lam + 1) \nl (\siro \batu \siro) + \lam \nl (\siro \siro \siro) + \nl (\siro) = 0.
\end{eqnarray}
\end{cor}
The detailed discussion concerning results in this section can be seen in Konno \cite{k1994,k1997}. If we regard $\lam, \nl (\siro), \nl (\siro \siro), \nl (\siro \siro \siro), \ldots $ as variables, then the left hand sides of the correlation identities by Theorem {\rmfamily \ref{thm1}} are polynomials of degree at most two. In the next section, we give a new procedure for getting a series of approximations for extinction probabilities based on the Gr\"obner basis by using Corollary {\rmfamily \ref{cor2}}. As for the Gr\"obner basis, see \cite{clo1992}, for example.

\section{Our results}
Put $x = \nl (\siro), \> y = \nl (\siro \siro), \> z = \nl (\siro \siro \siro), \> w = \nl (\siro \batu \siro), \> s = \nl (\siro \siro \siro \siro), \> u = \nl (\siro \siro \batu \siro)$. Let $\prec$ denote the lexicographic order with $\lam \prec x \prec y \prec w \prec z \prec u \prec s.$ For $m=1,2,3$, let $I_m$ be the ideals of a polynomial ring $\RM [x_1, x_2, \ldots, x_{n(m)}]$ over $\RM$ as defined below. Here $x_1 = \lam, x_2 = x, x_3 = y, x_4 = z, x_5 = w, x_6 = s, x_7 = u$ and $n(1)=3, n(2)=4, n(3)=7.$

\subsection{First approximation}
We consider the following ideal based on Corollary {\rmfamily \ref{cor2}} (1):
\begin{eqnarray}
I_1 = \langle \> 2 \lam y - 2 \lam x - x + 1, \> y - x^2 \> \rangle \subset \RM [\lam,x,y].
\end{eqnarray}
Here $y - x^2$ corresponds to the first (or mean-field) approximation: $\nl^{(1)}(\siro \siro) = (\nl^{(1)}(\siro))^2$. Then 
\begin{eqnarray}
G_1 = \{ (x-1)(2 \lam x -1), \> y-x^2 \}
\end{eqnarray}
is the reduced Gr\"obner basis for $I_1$ with respect to $\prec$. Therefore the solution except a trivial one $x(=y)=1$ is $x=\nl^{(1)} (\siro)= 1/(2 \lam)$. Remark that the trivial solution means that the invariant measure is $\delta_0$. From this, we obtain the first approximation of the density of the particle, $\rho_{\lam} = E_{\nl}(\eta (x))$, as follows:
\begin{eqnarray}
\rho_{\lam} ^{(1)} = 1- \nl^{(1)} (\siro) = {2 \lam -1 \over 2 \lam},
\end{eqnarray}
for any $\lam \ge 1/2$. This result gives the first lower bound $\lam_c ^{(1)}$ of the critical value $\lam_c$ of the one-dimensional contact process, that is, $\lam_c ^{(1)} =1/2 \le \lam_c.$ However it should be noted that the inequality is not proved in our approach. The estimated value of $\lam_c$ is about 1.649. 

\subsection{Second approximation}
Consider the following ideal based on Corollary {\rmfamily \ref{cor2}} (1) and (2):
\begin{eqnarray*}
I_2 = \langle \> 2 \lam y - 2 \lam x - x + 1, \> \lam z - \lam y - y + x, \> xz - y^2 \> \rangle \subset \RM [\lam,x,y,z].
\end{eqnarray*}
Here $xz - y^2$ corresponds to the second (or pair) approximation: $\nl^{(2)}(\siro) \nl^{(2)}(\siro \siro \siro) = (\nl^{(2)}(\siro \siro))^2$. Then 
\begin{eqnarray*}
&&
G_2 = \{ (x-1)((2 \lam -1) x -1), \> 1+2 \lam(y-x)-x, \\
&&
\qquad \qquad \qquad \qquad 
\> -y -yx + 2x^2, -z -y(2+y)+4x^2 \}
\end{eqnarray*}
is the reduced Gr\"obner basis for $I_2$ with respect to $\prec$. Therefore the solution except a trivial one $x(=y=z)=1$ is $x = \nl^{(2)} (\siro) = 1/(2 \lam-1)$. As in a similar way of the first approxaimation, we get the second approximation of the density of the particle:
\begin{eqnarray*}
\rho_{\lam} ^{(2)} = {2 (\lam -1) \over 2 \lam -1},
\end{eqnarray*}
for any $\lam \ge 1$. This result implies the second lower bound $\lam_c ^{(2)} =1$. We should remark that if we take 
\begin{eqnarray*}
I_2 '= \langle \> 2 \lam y - 2 \lam x - x + 1, \> \lam z - \lam y - y + x, \> y-x^2, \> z - x^3 \> \rangle \subset \RM [\lam,x,y,z],
\end{eqnarray*}
then we have
\begin{eqnarray*}
G_2 '= \{ z-1, \> y-1, \> x-1 \}
\end{eqnarray*}
is the reduced Gr\"obner basis for $I_2 '$ with respect to $\prec$. Here $y - x^2$ and $z - x^3$ correspond to an approximation: $\nl^{(2')}(\siro \siro) = (\nl^{(2')}(\siro))^2$ and $\nl^{(2')}(\siro \siro \siro) = (\nl^{(2')}(\siro))^3$, respectively. Then we have only trivial solution: $x=y=z=1.$ 

\subsection{Third approximation}
Consider the following ideal based on Corollary {\rmfamily \ref{cor2}} (1)--(4):
\begin{eqnarray*}
&&
I_3 = \langle \> 2 \lam y - 2 \lam x - x + 1, \> \lam z - \lam y - y + x, 
\\
&&
\qquad \qquad 
2 \lam s + w - (2 \lam + 3) z + 2 y, 
\> \lam u - (2 \lam + 1) w + \lam z + x, 
\\
&&
\qquad \qquad \qquad \qquad \qquad \qquad
y s - z^2, \> xu - yw \> \rangle \subset \RM [\lam,x,y,z,w,s,u].
\end{eqnarray*}
Here $ys - z^2$ and $xu - yw$ correspond to the third approximation: $\nl^{(3)}(\siro \siro) \nl^{(3)}(\siro \siro \siro \siro) = (\nl^{(3)}(\siro \siro \siro))^2$ and $\nl^{(3)}(\siro) \nl^{(3)}(\siro \siro \batu \siro) = \nl^{(3)}(\siro \siro) \nl^{(3)}(\siro \batu \siro)$, respectively. Then 
\begin{eqnarray*}
G_3 = \{ (x-1)((12 \lam^3 - 5 \lam -1) x^2 - 2 \lam (2 \lam +3)x - \lam +1), \> \ldots \}
\end{eqnarray*}
is the reduced Gr\"obner basis for $I_3$ with respect to $\prec$. Therefore the solution except a trivial one $x=1$ is $x= \nl^{(3)} (\siro) = (\lam (2 \lam +3) + \sqrt{D})/(12 \lam^3 - 5 \lam -1)$, where $D=16 \lam^4 + 4 \lam^2 + 4 \lam +1$. Then we obtain the third approximation of the density of the particle:
\begin{eqnarray}
\rho_{\lam} ^{(3)} = {4 \lam (3 \lam^2 - \lam -3) \over 12 \lam^3 - 2 \lam^2 -8 \lam -1 + \sqrt{D}},
\end{eqnarray}
for any $\lam \ge (1 + \sqrt{37})/6$. This result corresponds to the third lower bound $\lam_c ^{(3)} =(1 + \sqrt{37})/6 \approx 1.180$.

\section{Summary}
We obtain the first, second, and third approximations for the extinction probability, the density of the particle, and the lower bound of the one-dimensional contact process by using the Gr\"obner basis with respect to a suitable term order. These results coincide with results given by the Harris lemma (more precisely, the Katori-Konno method, see \cite{k1997}) or the BFKL inequality \cite{bfkl1997} (see also \cite{k1997}). As we saw, the generators of $I_m$ in Section 2 have degree at most two in $x_1,x_2, \ldots$, such as $2 \lam y - 2 \lam x - x + 1, \> y s - z^2$ in the case of $I_3$. We expect that this property will lead to get the higher order approximations of the process (and other interacting particle systems having a similar property) effectively. 
\par
\
\par\noindent
{\bf Acknowledgment.} The author thanks Takeshi Kajiwara for valuable discussions and comments.

\end{document}